\documentclass[a4paper,twoside]{article}

\usepackage{color}
\usepackage{authblk}
\usepackage{makeidx}         
\makeindex   
\usepackage{graphicx}        
\usepackage[all]{xy}
\usepackage{amssymb}
\usepackage{amsmath}
\usepackage{amsthm}
\usepackage{hhline}
\usepackage{longtable}
\usepackage[a4paper, left=2cm, right=2cm, top=2cm,bottom=2cm]{geometry}
\usepackage{hyperref}
\usepackage[draft=true]{minted}
\usepackage[backend=bibtex, bibencoding=utf8, doi=false, url=false, isbn=false]{biblatex}
\theoremstyle{definition}
\newtheorem{theorem}{Theorem}[section]
\newtheorem{definition}[theorem]{Definition}

\usepackage[all]{xy}
\usepackage{amssymb}
\usepackage{amsmath}

\newcommand{\NN}{\mathbb{N}}

\newcommand{\CC}{\mathbb C}

\newcommand{\QQ}{\mathbb Q}
\newcommand{\ZZ}{\mathbb Z}

\newcommand{\Aut}{\mathrm{Aut}}
\newcommand{\Triv}{\mathrm{Triv}}
\newcommand{\nef}{\mathrm{nef}}
\newcommand{\aut}{\mathrm{aut}}
\newcommand{\rk}{\mathrm{rank}}
\newcommand{\Div}{\mathrm{Div}}

\newcommand{\NS}{\mathrm{NS}}

\begin{document}

\title{Elliptic Fibrations on Vinberg's Most Algebraic K3 Surface}

\author[$\dagger$]{Simon Brandhorst}
\author[$\heartsuit$]{Matthias Zach}

\affil[$\dagger$]{
  Universit\"at des Saarlandes,
  66123 Saarbr\"ucken,
  Germany 
}
%
\affil[$\heartsuit$]{
  RPTU Kaiserslautern-Landau,
  67663 Kaiserslautern,
  Germany 
}


\maketitle

\abstract{We explain how to use the computer algebra system \verb|OSCAR| to find all elliptic fibrations (up to automorphism) on a given surface and
compute their Weierstrass models. This is illustrated for Vinberg's most algebraic K3 surface, the unique K3 surface of Picard rank $20$ and discriminant $3$.}

\section{Introduction}
\label{sec:introduction}

\verb|OSCAR|\cite{OSCAR} is a fairly new computer algebra system based on the \verb|Julia|\cite{julia} 
programming language. It incorporates and extends several established systems and libraries such as 
\verb|Singular|\cite{SINGULAR}, \verb|GAP|\cite{GAP4}, \verb|polymake|\cite{polymake}, 
\verb|Nemo|\cite{MR3703682}, and others. The joint power of these subsystems together with the 
versatility of the \verb|Julia| language allows for new efficient modeling of algebraic 
geometry way beyond the cases where an affine scheme is represented by a single Noetherian ring $R$ or 
a projective scheme by a graded unital algebra $S$. The framework for algebraic schemes 
in \verb|OSCAR| is currently under active development as part of the SFB-TRR 195 by the 
DFG (Deutsche Forschungsgemeinschaft). In this note we wish to illustrate the power it 
has already gained by studying elliptically fibered K3 surfaces and their automorphisms. 
The code used has been integrated in \verb|OSCAR|, version \verb|0.14.0-DEV|.

\medskip

A K3 surface over a field $k$ is a complete non-singular surface \footnote{separated, geometrically integral scheme of finite type over $k$ and dimension $2$} $X/k$ with
\[ h^1(X, \mathcal{O}_X)=0, \quad \omega_X \cong \mathcal{O}_X.\]
Examples are double covers of the projective plane \(\mathbb{P}^2\) branched over a smooth sextic, smooth quartics in $\mathbb{P}^3$,
smooth complete intersections of a quadric and a cubic in $\mathbb{P}^4$.
If $A$ is an abelian surface, e.g. the Jacobian of a genus $2$ curve, then the minimal model of $A/{\pm1}$ is a K3 surface, a so called Kummer K3 surface.

K3 surfaces have a central place in algebraic geometry. They can be seen as a generalization of elliptic curves and as the two-dimensional incarnation of Calabi-Yau and Hyperk\"ahler manifolds. As such they have been studied extensively over the past decades and their geometry is by now reasonably well understood, see \cite{Huybrechts2016}. For instance the Tate conjecture (see \cite{Charles2013,Maulik2014, MadapusiPera2015,Pera2016} for $k$ a finite field) and the cone conjecture \cite{Sterk1985, Maulik2018} are both known for K3 surfaces
yet wide open in general. For \textit{elliptically fibered}\footnote{See Definition \ref{def:GenusOneFibration} below} K3 surfaces, the situation is under even better control; for example the
Tate conjecture for elliptic K3 surfaces was known already in the '70s \cite{ArtinSwinnerton1973}.

A distinctive feature of K3 surfaces is that a single K3 surface $X$
may admit several elliptic fibrations. Such a fibration $\pi\colon X
\to \mathbb{P}^1$ is determined by the class of a fiber of $\pi$ in the
N\'eron-Severi lattice $\NS(X)$. We call it an \emph{elliptic class} \footnote{See Section \ref{subsec:fibrations_and_automorphisms} for a definition}.
Shioda and Kuwata \cite{KuwataShioda2008} posed the following problems:\\
Given an elliptic class
$f \in \NS(X)$ in the N\'eron-Severi group of $X$
determine
\begin{enumerate}
 \item[(1)] the elliptic parameter $t$ for $f$;
 \item[(2)] the defining equation of the elliptic curve $E/k(t)$;
 \item[(3)] the Mordell-Weil lattice of $E/k(t)$ and its generators;
 \item[(4)] all the elliptic fibrations $\pi': X \to \mathbb{P}^1$ up to isomorphism.
\end{enumerate}
These problems have been solved to varying degree on several (families of) K3 surfaces.
The automorphism group of a Kummer K3 surface $X$ of the Jacobian of a genus two curve $C$ (whose Jacobian $J(C)$ has no extra endomorphisms)
was computed by Kondo in \cite{Kondo1998}. This was used by Kumar in \cite{Kumar2014} to compute
representatives for the set of elliptic (and genus one) fibrations on $X$ modulo $\Aut(X)$
and their elliptic parameters. There are 25 elliptic fibrations on a very general Jacobian Kummer K3 $X$ up to automorphisms.
Similar
classifications on other (families of) K3 surfaces have been worked out by
Elkies and Sch\"utt \cite{ElkiesSchuett2015}, Nishiyama \cite{Nishiyama1996}, Kuwata and Shioda \cite{KuwataShioda2008} largely by hand. In \cite{FestiVeniani2022}, Festi and Veniani give a formula for the number of elliptic fibrations up to automorphisms.
The first named author and Elkies used the techniques to find equations for Lehmer's automorphism of minimal topological entropy on a K3 surface.

We demonstrate how to solve Shioda's and Kuwata's problems using Oscar by considering Vinberg's 'most algebraic' K3 surface $Y_1/\mathbb{C}$. The proofs of the algorithms involved and the details of the implementation shall appear elsewhere.


\section{The most algebraic K3 surface}
Vinberg's most algebraic K3 surface $Y_1$ is the unique complex K3 surface of Picard number $20$ and discriminant $3$ \cite{vinberg1983}. In what follows we will illustrate (some part of) the theory of K3 surfaces and their elliptic fibrations by considering $Y_1$ as a running example.

\subsection{Elliptic Fibrations}
In this section we fix our notation on elliptic surfaces and recall some definitions.
For a detailed account of the theory of elliptic surfaces we refer to \cite{Miranda1989, SchuettShioda2019}. For simplicity, we assume throughout that the characteristic of the base field $k$ is not $2$ or $3$.
\begin{definition}
  \label{def:GenusOneFibration}
A genus-one fibration on a surface $X$ is a surjective map
\(\pi \colon X \to C\)
to a smooth curve $C$ such that the generic fiber is a smooth curve of genus one.
If $\pi$ has a section $O$, we call the triple $(X,\pi,O)$ an elliptic fibration and $O$ the zero section.
\end{definition}
Any choice of a zero section turns the generic fiber of $\pi$ into an elliptic curve $E$ over $K=k(C)$, the function field of $C$. If $\pi$ and $O$ are clear from the context, they are omitted from the notation and we call $X$ simply an elliptic surface.
For a K3 surface $X$ only the base $C = \mathbb{P}^1$ is possible because any regular one form on the base pulls back to a regular 1-form on $X$ and K3 surfaces have no regular $1$-forms.

We say that $(X,\pi)$ is \textit{relatively minimal}, if $X$ is smooth and the fibers of $\pi$ do not contain any $(-1)$-curves.
Since such curves can be contracted to a smooth point, an arbitrary genus-one fibration $(X', \pi')$ 
with smooth total space $X'$ can be blown down to a relatively minimal one. Conversely, given an elliptic curve $E/k(C)$ one can construct a relatively minimal elliptic fibration whose generic fiber is $E$.
This elliptic fibration is called the Kodaira-N\'eron model of $E$ and is unique up to isomorphisms \cite[Thm. 5.19]{SchuettShioda2019}. For $C = \mathbb{P}^1$ we will see later how to construct it in Oscar.

Let
\begin{equation}\label{eq:generic_fiber}
E_1/\CC(t): y^2 = x^3 +t^5(t-1)^2
\end{equation}
and \(\pi\colon Y_1 \to \mathbb{P}^1\) the
relatively minimal elliptic fibration defined by $E_1$.
Then we will see that \(Y_1\) is Vinberg's most algebraic K3 surface.

In Oscar we define $E_1$ as an elliptic curve over $\QQ(t)$.

\begin{minted}{jlcon}
julia> Qt, t = polynomial_ring(QQ, :t)
(Univariate polynomial ring in t over QQ, t)

julia> Qtf = fraction_field(Qt)
Fraction field
  of univariate polynomial ring in t over QQ

julia> E = EllipticCurve(Qtf, [0,0,0,0,t^5*(t-1)^2])
Elliptic curve with equation
y^2 = x^3 + t^7 - 2*t^6 + t^5
\end{minted}
The Kodaira-N\'eron Model $Y_1$ of $E_1$ has a contraction to a Weierstrass model $S_1$ given by the homogeneous equation $Y^2 = X^3 + t^5(t-s)^2s^5$ of degree $12$ in weighted projective space $\mathbb{P}(1,1,2,3)$ with coordinates $(s,t,X,Y)$. Later we construct $S$ in a certain projective bundle of $\mathbb{P}^1$. At this point we remark that $S_1$ has an $E_8$ singularity in the fibers over $[0:1]$ and $[1:0]$ each and an $A_2$ singularity in the fiber over $[1:1]$. The minimal resolution of these singularities is the contraction $Y_1\to S_1$. The exceptional divisors are contained in the respective fibers. This is how reducible fibers may appear. They will play an important role in the next section.
\noindent

\subsubsection{Fibers and the trivial lattice}\label{subsec:mwl}
Since any complex K3 surface $X$ is simply connected, the
notions of linear, algebraic and numerical equivalence agree
and yield $\operatorname{Pic}(X) \cong \operatorname{NS}(X)
\cong \operatorname{Num}(X)$ \cite{Huybrechts2016}. This is also true in positive
characteristic. Equipped with the intersection pairing
$\operatorname{NS}(X)$ is a
\textit{quadratic integer lattice}.
Its rank is called the Picard number of $X$.
Computing generators of $\operatorname{NS}(X)$ or even its rank is in general a hard problem --- just like computing rational points on elliptic curves.

On an elliptic fibration some part of $\operatorname{NS}(X)$ is readily computable,
the so called \emph{trivial lattice}
\(\operatorname{Triv}(X,\pi,O)\) spanned by the zero section and components of the
reducible fibers.
The classes $o,f \in \operatorname{NS}(X)$ of the zero section $O$ and any fiber $F$ of $\pi$ satisfy $f^2=0$, $o^2=-2$
and $f.o=1$, i,e. they span a so called
\textit{hyperbolic plane} in $\NS(X)$.
An elliptic fibration may have reducible fibers.
The reducible fibers of $\pi$ which are
disjoint from the zero section \(O\) decompose as a direct sum
of (negative definite) ADE-root lattices --- as predicted by Kodaira
classification of singular fibers \cite[Thm. 6.2]{Kodaira63} and Tate's algorithm \cite{Tate1975} to compute
them. Tate's algorithm implies that
we can compute the Kodaira types of the singular fibers by factoring the
$j$-invariant and the discriminant of $E$ and comparing this with the
table \cite[p.365]{silverman1994} (outside of characteristic $2$ and $3$).
In our case $j(E)=0$.
\begin{minted}{jlcon}
julia> factor(discriminant(E), Qt)
-432 * (t - 1)^4 * t^10
\end{minted}
\noindent
Since the base of the fibration is $\mathbb{P}^1$, we should view the
discriminant as a homogeneous polynomial of degree $12 l$ in
$[t:s]$ where $l>0$ comes from the construction of the Kodaira-N\'eron
model. For $X$ a K3 surface, $l=2$ and the discriminant has degree
$24$ which is -- and this is not a coincidence --
the topological Euler number of any K3 surface.

We are in the affine chart $s=1$ and therefore the homogeneous
polynomial for the discriminant is 
\[\Delta(E) = -432 (t-s)^4 t^{10} s^{10}.\]
We see that the fibration $\pi$ has two singular fibers of type \(II^*\) (an extended $E_8$ configuration) at $[0:1]$ and $[1:0]$, as well as a fiber of type $IV$ (forming an extended $A_2$ configuration) at $[1:1]$.

The components not meeting the zero section form an $E_8$, $E_8$, $A_2$ configuration.
Therefore the trivial lattice is
\(\Triv(\pi) \cong U \oplus E_8 \oplus E_8 \oplus A_2\)
which is of rank $20$ and determinant $-3$.
Since the maximum Picard number of a K3 surface in characteristic $0$ is $20$ and $-3$ is the minimal possible determinant, we conclude that \(\NS(X) = \Triv(\pi)\)
is the full N\'eron-Severi lattice. Therefore $Y_1$ is indeed Vinberg's most algebraic K3 surface \cite{vinberg1983}.

\subsection{Elliptic fibrations and the automorphism group}\label{subsec:fibrations_and_automorphisms}
A distinctive feature of K3 surfaces is that they may admit several elliptic fibrations.
They are in bijection with the extremal rays of the Nef cone.
\begin{theorem}\cite[\S 3 Theorem 1]{Pjatecki1971}\cite[Lemma 11.25]{SchuettShioda2019}
 Let $X$ be a K3 surface and $f \in \operatorname{NS}(X)$ a primitive and Nef divisor class with $f^2=0$. Then the map $\pi_{|f|}$ induced by the linear system $|f|$ is a genus $1$ fibration.
 The ADE-type of its singular fibers is given by the ADE-type of the root sublattice of $(f^\perp/ \ZZ f)$.
 If there is a class $e \in \NS(X)$ with $e.f = 1$, then the fibration has a section.
\end{theorem}
We call a primitive and Nef divisor class $f \in \NS(X)$ an \textit{elliptic class} and
$\pi^*(t/s) \in k(X)$ an elliptic parameter where $t/s$ is a generator of the function field of the base $\mathbb{P}^1$.
Let us point out that if $f\neq0$ is isotropic but not Nef, then it can be mapped to an elliptic class $f'$ by an element of the Weyl group. Therefore, by Mayer's theorem, any K3 surface of Picard number at least $5$ has a genus one fibration.

To classify all elliptic fibrations on $X$ up to automorphism, we need a description of its automorphism group.
Over $\CC$ the strong Torelli type theorem implies that $\Aut(X)$ can
be identified with the group of effective, integral Hodge isometries of
$H^2(X,\ZZ)$, i.e. isometries that preserve the intersection product,
the Nef cone $\nef(X)$ and the integral Hodge structure on $H^2(X,\ZZ)$. We refer to \cite[Section 7]{Huybrechts2016} for details.
More generally, for an arbitary field $k$ and any K3 surface $X/k$ the map
\begin{eqnarray}
 \Psi\colon \Aut(X) &\to &O(\NS(X), \nef(X))\cong O(\NS(X))/W(\NS(X)) \\
 \phi &\mapsto &(\phi^*)^{-1}
\end{eqnarray}
has finite kernel and cokernel. Here $W(\NS(X))$ denotes the Weyl group, that is the group generated by all reflections in $(-2)$-vectors of $\NS(X)$.
Denote by $\aut(X)$ the image of $\Psi$. For \(Y_1/\mathbb{C}\) we know (from the
K3-Torelli theorem) that $\Psi$ is an isomorphism.
To solve Shioda and Kuwata's problems, we have to compute $\aut(X)$, $\nef(X)$, and
\[
  \{f \in \NS(X) \, | \, f \mbox{ is elliptic }\}/\aut(X).
\]

Using an implementation in Oscar of Shimada's algorithm
\cite{Shimada2015}, which builds on ideas of Borcherds' and Vinberg,
we compute $\aut(Y_1)$ as well as a fundamental domain of its action on
the Nef cone.

First we set up the N\'eron-Severi lattice.
\begin{minted}{jlcon}
julia> R = rescale(root_lattice([(:E, 8), (:E, 8), (:A, 2)]), -1);

julia> U = integer_lattice(gram=ZZ[0 1; 1 -2]);

julia> NS, _ = direct_sum([U, R]);
\end{minted}
\noindent
Next we need to find an ample class, to identify the ample cone as a Weyl-chamber.
The ample class must intersect each curve (in particular the ones corresponding to the basis vectors) positively.
\begin{minted}{jlcon}
julia> e = matrix(ZZ,1,20,ones(Int,20));e[1,1]=51;

julia> ample = e*inv(gram_matrix(NS));
\end{minted}
\noindent
We confirm that the class is indeed in the interior of a Weyl-chamber,
by checking that it is not contained in any hyperplane orthogonal to a $(-2)$-class.
\begin{minted}{jlcon}
julia> minimum(rescale(orthogonal_submodule(NS, ample),-1))
4
\end{minted}
\noindent
We let Oscar compute generators for the automorphism group $\aut(Y_1)$.
\begin{minted}{jlcon}

julia> Xaut, Xchambers, Xrational_curves =  K3_surface_automorphism_group(NS, ample);
\end{minted}
\noindent
The return value \verb|Xchambers| is a list of polyhedral cones $C_1,\dots C_n$ which are contained in the Nef cone.
Let $\Aut(C_i):=\Aut(X)_{C_i}$ be the stabilizer of $C_i$. This stabilizer acts on the points of $C_i$. Let $P_i\subseteq C_i$ be its fundamental domain.
A fundamental domain of the action of $\Aut(X)$ on the Nef cone is the union of the $P_i$.
In our case $n=1$.
\begin{minted}{jlcon}
julia> C = Xchambers[1]
Chamber  in dimension 20 with 36 walls
\end{minted}
\noindent
The walls $W=r^\perp$ of $C$ with $r^2=-2$ and $r \in \NS(X)$ are classes
of smooth rational curves. The rays $f$ of $C$ with $f^2=0$ give elliptic
fibrations.  We now compute the rays of $C$ modulo $\Aut(C)$.
\begin{minted}{jlcon}
julia> Xelliptic_classes = [f for f in rays(C) if 0 == f*gram_matrix(NS)*transpose(f)];

julia> Xelliptic_classes_orb = orbits(gset(matrix_group(aut(C)),(x,g)->x*matrix(g), Xelliptic_classes));

julia> length(Xelliptic_classes_orb)
6
\end{minted}
\noindent
Thus modulo $\Aut(Y_1/\mathbb{C})$ there are exactly $6$ elliptic fibrations on $Y_1/\mathbb{C}$, this recovers a result of Nishiayama, Watari and Braun (see the introduction \cite{FestiVeniani2022})
We will see later they are all defined over $\mathbb{Q}$ -- because so is the full N\'eron-Severi-lattice of $Y_1$.

\medskip

In what follows we aim to compute different Weierstrass-models of
Vinberg's K3 surface which realize the different elliptic fibrations. 
Moreover, we wish to identify their minimal desingularizations by explicit rational maps.
It turns out that such a rational map can be arbitrarily complicated, meaning that only writing down the required rational functions to describe them might fill several pages.
Therefore, it is preferable to work with a class of simple, better controlled isomorphisms and
factor the more complicated ones in terms of these. A first bound on the complexity of a rational map between two elliptic fibrations is the intersection number
of their associated elliptic classes:

\begin{definition}
 Let $p \in \NN$ be prime. Two elliptic classes $f_1,f_2 \in \NS(X)$ are called $p$-neighbors if $f_1.f_2=p$.
\end{definition}
In fact one can show that the corresponding frame lattices of $f_1$ and $f_2$ can be seen as $p$-neighbors in the sense of Kneser if $\det(\NS(X))$ is coprime to $p$ \cite{BrandhorstElkies2023}.

The simplest non-trivial case is when two elliptic classes are $2$-neighbors.
Therefore, to compute all the sought for Weierstrass models and the corresponding elliptic parameters, we first
single out $6$ representatives $f_1,\dots f_6$, which
we can reach from our initial fibration by
$2$-neighbor steps.
\begin{minted}{jlcon}
julia> f1 = identity_matrix(ZZ, 20)[1,:];

julia> f1_2neighbors = [f for f in Xelliptic_classes if 2 == (f1*gram_matrix(NS)*transpose(f))[1,1]];

julia> f2, f3, _ = f1_2neighbors;

julia> f2_2neighbors = [f for f in Xelliptic_classes if 2 == (f2*gram_matrix(NS)*transpose(f))[1,1]];

julia> fibers = [f1,f2,f3];

julia> todo = [o for o in Xelliptic_classes_orb if !any(f in o for f in fibers)];

julia> while length(todo) > 0
         o = popfirst!(todo)
         f = findfirst(x-> x in o, f2_2neighbors)
         push!(fibers, f2_2neighbors[f])
       end

julia> fibers = [vec(collect(i*basis_matrix(NS))) for i in fibers]
6-element Vector{Vector{QQFieldElem}}:
 [1, 0, 0, 0, 0, 0, 0, 0, 0, 0, 0, 0, 0, 0, 0, 0, 0, 0, 0, 0]
 [4, 2, -4, -7, -10, -8, -6, -4, -2, -5, -2, -4, -6, -5, -4, -3, -2, -3, -1, -1]
 [4, 2, -4, -7, -10, -8, -6, -4, -2, -5, -4, -7, -10, -8, -6, -4, -2, -5, 0, 0]
 [6, 3, -5, -10, -15, -12, -9, -6, -3, -8, -5, -10, -15, -12, -9, -6, -3, -8, -1, -1]
 [10, 5, -8, -16, -24, -20, -15, -10, -5, -12, -7, -14, -21, -17, -13, -9, -5, -11, -3, -2]
 [6, 3, -4, -8, -12, -10, -8, -6, -3, -6, -4, -8, -12, -10, -8, -6, -3, -6, -2, -1]
\end{minted}
\noindent
Now we compute their Mordell-Weil rank, the torsion subgroup of the Mordell-Weil-group and the
ADE-type of their singular fibers. This recovers a result of Nishiyama \cite[Table 1.1]{Nishiyama1996}.
\begin{minted}{jlcon}
julia> [fibration_type(NS, f) for f in fibers]
6-element Vector{Tuple{Int64, GrpAbFinGen, Vector{Tuple{Symbol, Int64}}}}:
 (0, GrpAb: Z/1, [(:A, 2), (:E, 8), (:E, 8)])
 (1, GrpAb: Z/2, [(:D, 10), (:E, 7)])
 (0, GrpAb: Z/2, [(:A, 2), (:D, 16)])
 (1, GrpAb: Z/3, [(:A, 17)])
 (0, GrpAb: Z/4, [(:A, 11), (:D, 7)])
 (0, GrpAb: Z/3, [(:E, 6), (:E, 6), (:E, 6)])
\end{minted}
\noindent
We compute their mutual intersections $(f_i.f_j)_{1\leq i,j \leq 6}$
\begin{minted}{jlcon}
julia> F = transpose(matrix((reduce(hcat, fibers)))); F * gram_matrix(NS) * transpose(F)
[0   2   2   3   5   3]
[2   0   2   2   2   2]
[2   2   0   2   5   4]
[3   2   2   0   2   3]
[5   2   5   2   0   2]
[3   2   4   3   2   0]
\end{minted}
\noindent and see that $f_1.f_2=f_1.f_3=2$ and
$f_2.f_i = 2$ for $i\neq 2$. Therefore, we can reach $f_2,f_3$ by a single $2$-neighbor step from $f_1$ and $f_4,f_5,f_6$ by another one from $f_2$.

\subsection{The (relatively) minimal model}
To solve the next problem on Shimada's list, we return from lattices back to algebraic geometry.

In Oscar we can construct $\pi: Y_1 \to \mathbb{P}^1$ as a covered scheme which is represented internally by patches and glueing data, just like in Harshorne's book Algebraic Geometry \cite{Hartshorne1977}. Recall that $E_1/\QQ(t)$ is the generic fiber of $\pi_1$.
The (singular) Weierstrass model $S_1$ of $(Y_1,\pi_1)$ is the zero locus of a section of the projective bundle $\mathbb{P}(\mathcal{O}_{\mathbb{P}^1}(-2l) \oplus \mathcal{O}_{\mathbb{P}^1}(-3l) \oplus \mathcal{O}_{\mathbb{P}^1}(1))$ where $l=2$ since $Y_1$ is a K3 surface \cite[\S 5.13]{SchuettShioda2019}.
\begin{minted}{jlcon}
julia> Y1 = elliptic_surface(E, 2)
Elliptic surface
  over rational field
with generic fiber
  -x^3 + y^2 - t^7 + 2*t^6 - t^5

julia> S = weierstrass_model(Y1)[1]
Scheme
  over rational field
with default covering
  described by patches
    1: spec of quotient of multivariate polynomial ring
    2: spec of quotient of multivariate polynomial ring
    3: spec of quotient of multivariate polynomial ring
    4: spec of quotient of multivariate polynomial ring
    5: spec of quotient of multivariate polynomial ring
    6: spec of quotient of multivariate polynomial ring
  in the coordinate(s)
    1: [(x//z), (y//z), t]
    2: [(z//x), (y//x), t]
    3: [(z//y), (x//y), t]
    4: [(x//z), (y//z), s]
    5: [(z//x), (y//x), s]
    6: [(z//y), (x//y), s]
\end{minted}
\noindent
We compute a resolution of singularities $\psi \colon Y_1\to S_1$ by performing a series of blowups in one singular point at a time. This is fully automated and leads to a total of 17 blowups.
\begin{minted}{jlcon}
julia> piS = weierstrass_contraction(Y1)
(Scheme over QQ covered with 46 patches, Composition of Scheme over QQ covered with 46 patches -> Scheme over QQ covered with 44 patches -> Scheme over QQ covered with 42 patches -> Scheme over QQ covered with 38 patches -> Scheme over QQ covered with 36 patches -> Scheme over QQ covered with 34 patches -> Scheme over QQ covered with 32 patches -> Scheme over QQ covered with 30 patches -> Scheme over QQ covered with 28 patches -> Scheme over QQ covered with 26 patches -> Scheme over QQ covered with 24 patches -> Scheme over QQ covered with 20 patches -> Scheme over QQ covered with 18 patches -> Scheme over QQ covered with 16 patches -> Scheme over QQ covered with 14 patches -> Scheme over QQ covered with 12 patches -> Scheme over QQ covered with 10 patches -> Scheme over QQ covered with 6 patches)
\end{minted}
\noindent
Note that this is less elegant than Tate's algorithm but it works for any surface with isolated ADE singularities.
We use the resolution to compute the exceptional divisors, i.e. the trivial lattice,
\begin{minted}{jlcon}
julia> basisNSY1, gramTriv = trivial_lattice(Y1);
\end{minted}
\noindent as well as the singular fibers.
\begin{minted}{jlcon}
julia> [(i[1],i[2]) for i in reducible_fibers(Y1)]
3-element Vector{Tuple{Vector{QQFieldElem}, Tuple{Symbol, Int64}}}:
 ([1, 1], (:A, 2))
 ([0, 1], (:E, 8))
 ([1, 0], (:E, 8))
\end{minted}
Since the computation of the full N\'eron-Severi lattice is a difficult problem, we work with a sublattice called the algebraic lattice in Oscar. It is spanned by fiber components and the sections supplied to the constructor.
Recall that for Vinberg's K3 surface $Y_1$ the trivial lattice agrees with $\NS(Y_1)$, hence so does the algebraic lattice.
\begin{minted}{jlcon}
julia> basisNSY1, _, NSY1 = algebraic_lattice(Y1);
\end{minted}
\noindent
The first return value consists of the chosen basis of the algebraic lattice. For brevity we only show the first 5 elements.
\begin{minted}{jlcon}
julia> basisNSY1
20-element Vector{Any}:
 Fiber over (2, 1)
 section: (0 : 1 : 0)
 component A2_1 of fiber over (1, 1)
 component A2_2 of fiber over (1, 1)
 component E8_1 of fiber over (0, 1)
\end{minted}
\noindent
They are represented as effective Weil divisors in Oscar.
\begin{minted}{jlcon}
julia> basisNSY1[1]
Effective weil divisor Fiber over (2, 1)
  on scheme over QQ covered with 46 patches
with coefficients in integer Ring
given as the formal sum of
  1 * sheaf of ideals
\end{minted}

\subsection{Divisors and linear systems on elliptic K3 surfaces}
Our goal is to compute the linear systems $|f_i|$, $i\in \{2,3,4,6\}$. To this end we need a better control on the Picard group. Here the elliptic fibration comes to help.

Recall that $K = k(C)$ is the function field of the base curve of an elliptic fibration $\pi\colon X \to C$. The $K$-rational points $E(K)$ of $E$ form a finitely generated abelian group group, which is called the Mordell-Weil-group of $E$.
A section of the fibration $\pi$ is a map $s\colon C \to X$ with $\pi \circ s = id_C$. We identify $s$ with its image $s(C)$ and call it a section as well.
The generic point of a section is a $K$-rational point of $E$. Conversely, the closure of a $K$-rational point of $E$ in $X$ is a section. This gives a bijection between the set of $K$-rational points $E(K)$ and the set of sections of $\pi$. We use it to give the set of sections the structure of an abelian group which we denote by \(\operatorname{Mwl}( \pi) \) and its addition by $\oplus$.

We call a prime divisor \textit{vertical} if its image under $\pi$ is a point and otherwise \textit{horizontal}. Any fibration structures the N\'eron-Severi-lattice into a vertical part and a horizontal part. The advantage of an
\emph{elliptic} fibration is that a horizontal divisor is linearly equivalent to a sum of sections (modulo vertical divisors).
\begin{theorem}\cite[Theorem 1.3]{shioda1990}
 Let $\pi\colon X \to C$ be an elliptic fibration. Then there is an isomorphism of groups
 \[\operatorname{Mwl}(\pi) \cong \operatorname{NS}(X)/\operatorname{Triv}(X,\pi)\]
 and an isometry of lattices
 \[\operatorname{Mwl}(\pi)_{free} \cong \operatorname{NS}(X)/\overline{\operatorname{Triv}(X, \pi)}.\]
\end{theorem}
Here $\overline{\operatorname{Triv}(X, \pi)}$ denotes the primitive closure of the trivial lattice in 
$\operatorname{NS}(X)$.
This theorem results in the so called \textit{Shioda Tate formula} which relates the Picard number of $X$ with the Mordell-Weil rank of the generic fiber of $\pi$.
\[ \rho(X) = 2 + \sum_{\nu \in \mathbb{P}^1} (n_\nu-1) + \rk \operatorname{MW}(\pi)\]
where $n_\nu$ is the number of irreducible components of the fiber $\pi^{-1}(\nu)$.
In our case we find that $E/\QQ(t)$ has zero Mordell-Weil rank.

On the side of divisors the theorem tells us the following (see \cite[Lemma 5.1]{shioda1990}:
Let $A \in \operatorname{Div}(X)$ be an effective divisor, then there is a linear equivalence
\[
  A \sim V + (n-1)O + P
\]
where $V$ is a vertical divisor, i.e. a sum of fiber components, and $n = A.F$. Thus when we want to compute linear systems on an elliptic surface,
it suffices to consider the linear systems of this form.

We wish to find $F_i \in \Div(X)$ effective with $[F_i]= f_i \in \NS(X)$. We start with $F_2$.
So first we have to identify the abstract $\NS(X)$ that we computed with the actual geometric one
on $X$ spanned by divisors.
\begin{minted}{jlcon}
julia> b, I = Oscar._is_equal_up_to_permutation_with_permutation(gram_matrix(NS), gram_matrix(NSY1))
(true, [1, 2, 19, 20, 3, 4, 5, 6, 7, 8, 9, 10, 11, 12, 13, 14, 15, 16, 17, 18])

julia> @assert gram_matrix(NSY1) == gram_matrix(NS)[I,I]
\end{minted}
\noindent
This identification is just a permutation of the basis.

Since the divisor $V$ is vertical, it is dominated by some divisor of fibers $V\leq V' = n_1 V_1 + \dots n_k V_k$ where $V_i=\pi^*(p_i)$ is a fiber over a point $p_i \in \mathbb{P}^1$ for $i \in \{1,\dots k\}$.
Note that in our case $n=2$ because we work with two neighbors.
Oscar can compute $V'$ and $V$ for us. We display $V$.
\begin{minted}{jlcon}
julia> Oscar.horizontal_decomposition(Y1, fibers[2][I])[2]
Effective weil divisor
  on scheme over QQ covered with 46 patches
with coefficients in integer Ring
given as the formal sum of
  2 * component E8_7 of fiber over (0, 1)
  2 * section: (0 : 1 : 0)
  1 * component A2_0 of fiber over (1, 1)
  1 * component E8_0 of fiber over (1, 0)
  2 * component E8_4 of fiber over (0, 1)
  2 * component E8_3 of fiber over (0, 1)
  1 * component E8_2 of fiber over (0, 1)
  2 * component E8_0 of fiber over (0, 1)
  2 * component E8_6 of fiber over (0, 1)
  2 * component E8_5 of fiber over (0, 1)
  1 * component E8_8 of fiber over (0, 1)
\end{minted}
\noindent
The larger linear system $|V' + O + P|$ is now of a simple form and computed in \cite{BrandhorstElkies2023}. Then Oscar can figure out the linear equations cutting out the subsystem $|V+O+P|$.
Since $F_2$ is an elliptic class, this linear system is of (projective) dimension one generated by two sections $v,w \in k(Y_1)$. Their quotient $u=v/w \in k(Y_1)$ is an elliptic parameter. The fibration $\pi_2 \colon Y_1 \to \mathbb{P}^1$ defined by $F_2$ is thus given by $P\mapsto (1:u(P))$
\begin{minted}{jlcon}
julia> representative(elliptic_parameter(Y1, fibers[2][I]))
(x//z)//(t^3 - t^2)
\end{minted}
with $u=x/(t^3-t^2)$ in Weierstrass coordinates.

\subsection{A 2-neighbor step}
To compute the fibrations of $F_i$ for $i=4,5,6$, we need to pass
to a Weierstrass model of the fibration $F_2$.
Using the elliptic parameter $u_i=u_i(x,y,t)$ Oscar can find a suitable
change of coordinates such that the generic fiber of $F_i$ is given as
double cover of $\mathbb{P}^1_{k(u_i)}$ ramified over four points.
An affine chart of this double cover is given by $g_i=y_i^2-q(x_i) \in k(u_i)[x_i,y_i]$ with $q_i(x_i)$ a separable polynomial degree $3$ or $4$.
It describes a curve of genus one over $k(u_i)$, but not yet an elliptic curve, because
in general it may not have a rational point. The following command
produces a pair \verb|(g, phi1)| consisting of $g_2$, together with a morphism $\varphi_1 \colon k[x, y, t] \to
(k(u_2))(x_2, y_2)$ taking the Weierstrass equation for the original $Y_1$
to an element with $g_2$ as the relevant factor of the numerator.
\begin{minted}{jlcon}
julia> g, phi1 = two_neighbor_step(Y1, fibers[2][I]); g
-t^3*x^3 + t^3*x^2 - x + y^2
\end{minted}
\noindent To compute a Weierstrass form we choose the rational point $P$, given by $x_2=y_2=0$,
of the generic fiber and move it to infinity. This produces an elliptic 
surface $Y_2$ and another homomorphism $\varphi_2 : (k(t_2))(x_2, y_2) \to k(Y_2)$
encoding the necessary change of coordinates:
\begin{minted}{jlcon}
julia> kt2 = base_ring(parent(g)); P = kt2.([0,0]);

julia> Y2, phi2 = elliptic_surface(g, P); Y2
Elliptic surface
  over rational field
with generic fiber
  -x^3 + t^3*x^2 - t^3*x + y^2
and relatively minimal model
  scheme over QQ covered with 47 patches
\end{minted}
From our computations in Section \ref{subsec:fibrations_and_automorphisms}
we already know that the corresponding elliptic fibration $\pi_2$ has a
$D_{10}$ and an $E_7$ fiber, a $2$-torsion section and Mordell-Weil rank
one. To find a section of infinite order, we note that $f_2$ intersects
some fiber components of $\pi_1$ with multiplicity one. Therefore they
are sections of $\pi_2$. We inform Oscar about the section.
\begin{minted}{jlcon}
julia> E2 = generic_fiber(Y2); t2 = gen(kt2);

julia> P2 = E2([t2^3, t2^3]); set_mordell_weil_basis!(Y2, [P2]);
\end{minted}
\noindent The upshot is that we have two elliptic surfaces $\pi_i\colon  Y_i \to \mathbb{P}^1$, $i=1,2$, and a birational map $\phi\colon  Y_2 \dashrightarrow Y_1$ which is given in Weierstrass coordinates of $Y_1$ and $Y_2$ by
\[(x,y,t) \mapsto \left( \frac{-t x + t}{x^3}, \frac{xy - y}{x^5}, \frac{1}{x}\right).\]
With Oscar this goes as follows.
\begin{minted}{jlcon}
julia> U2 = weierstrass_chart(Y2); U1 = weierstrass_chart(Y1);

julia> imgs = phi2.(phi1.(ambient_coordinates(U1))) # k(Y1) -> k(Y2)
3-element Vector{AbstractAlgebra.Generic.Frac{QQMPolyRingElem}}:
 (-(x//z)*t + t)//(x//z)^3
 ((x//z)*(y//z) - (y//z))//(x//z)^5
 1//(x//z)

julia> phi_rat = Oscar.MorphismFromRationalFunctions(Y2, Y1, U2, U1, imgs);
\end{minted}
Since the surfaces in question have Nef canonical class, we know from theory that the birational map is an isomorphism. Let us tell Oscar that $\phi$ is a morphism.
\begin{minted}{jlcon}
julia> set_attribute!(phi_rat, :is_isomorphism=>true);
\end{minted}
However actually computing $\phi$ as a morphism on some given chart can be
prohibitively expensive. This poses some serious computational challenges.

Our next task is to compute the matrix representing the pullback
$\phi^*\colon \NS(Y_1) \to \NS(Y_2)$ this is done by pulling
back the equations on carefully chosen charts. Then we use the intersection pairing.
\begin{minted}{jlcon}
julia> pullbackDivY1 = [pullback(phi_rat, D) for D in basisNSY1];

julia> B = [basis_representation(Y2, D) for D in pullbackDivY1];
\end{minted}
\noindent
At this point \verb|B| is a list of rational vectors. We turn it into 
a matrix and perform some sanity checks.
\begin{minted}{jlcon}
julia> B = matrix(QQ, 20, 20, reduce(vcat, B)); NSY2 = algebraic_lattice(Y2)[3];

julia> NSY1inY2 = lattice(ambient_space(NSY2),B);

julia> @assert NSY1inY2 == NSY2 && gram_matrix(NSY1inY2) == gram_matrix(NSY1)
\end{minted}
\noindent
Finally, we have the basis representation of the $f_i$ in $\NS(Y_2)$.
\begin{minted}{jlcon}
julia> fibers_in_Y2 = [f[I]*B for f in fibers]
6-element Vector{Vector{QQFieldElem}}:
 [4, 2, 0, 0, 0, 0, 0, 0, 0, -4, -4, -8, -7, -6, -5, -4, -3, -2, -1, 0]
 [1, 0, 0, 0, 0, 0, 0, 0, 0, 0, 0, 0, 0, 0, 0, 0, 0, 0, 0, 0]
 [5, 2, -2, -3, -4, -3, -2, -1, -2, -5, -4, -8, -7, -6, -5, -4, -3, -2, -1, 0]
 [4, 2, -2, -4, -6, -9//2, -3, -3//2, -7//2, -3, -5//2, -5, -9//2, -4, -7//2, -3, -5//2, -2, -3//2, 0]
 [2, 1, -1, -2, -3, -2, -1, 0, -2, -1, -1, -2, -2, -2, -2, -2, -2, -1, 0, 1]
 [2, 1, 0, 0, 0, 0, 0, 0, 0, -2, -2, -4, -4, -4, -4, -3, -2, -1, 0, 1]
\end{minted}

\subsection{The remaining fibrations}
With computations similar to the above we can now compute the elliptic parameters as well as the generic fibers for
$\pi_i \colon Y_i \to \mathbb{P}^1$ for $i=3,\dots 6$.

First we compute $\pi_3$ starting from $Y_1$.
\begin{minted}{jlcon}
julia> f3 = fibers[3][I]; representative(elliptic_parameter(Y1, f3))
(x//z)//t^2

julia> g3a, phi3a = two_neighbor_step(Y1, f3); g3a
-x^3 + (-t^3 + 2)*x^2 - x + y^2
\end{minted}
\noindent
Since we have the pullback action $\phi^*$ we are in a position to compute $\pi_i$ for $i=4,5,6$ from $Y_2$.
\begin{minted}{jlcon}
julia> [representative(elliptic_parameter(Y2, f)) for f in fibers_in_Y2[4:6]]
3-element Vector{AbstractAlgebra.Generic.Frac{QQMPolyRingElem}}:
 (y//z)//((x//z)*t)
 ((y//z) + t^3)//((x//z)*t - t^4)
 ((y//z) + t^3)//((x//z) - t^3)

julia> g,_ = two_neighbor_step(Y2, fibers_in_Y2[4]);g
-1//4*x2^4 - 1//2*t2^2*x2^3 - 1//4*t2^4*x2^2 + x2 + y2^2

julia> R = parent(g); K_t = base_ring(R);

julia> (x,y) = gens(R); P = K_t.([0,0]); # rational point

julia> g, _ = transform_to_weierstrass(g, x, y, P);

julia> E4 = elliptic_curve(g, x, y)
Elliptic curve with equation
y^2 = x^3 + 1//4*t2^4*x^2 - 1//2*t2^2*x + 1//4

julia> R = parent(g); K_t = base_ring(R);

julia> (x,y) = gens(R); P = K_t.([0,0]); # rational point

julia> g, _ = transform_to_weierstrass(g, x, y, P);

julia> E5 = elliptic_curve(g, x, y)
Elliptic curve with equation
y^2 = x^3 + (1//4*t2^4 - 2*t2)*x^2 + t2^2*x

julia> g,_ = two_neighbor_step(Y2, fibers_in_Y2[6]);g
(t2^2 + 2*t2 + 1)*x2^3 + y2^2 - 1//4*t2^4

julia> R = parent(g); K_t = base_ring(R); t = gen(K_t);

julia> (x,y) = gens(R); P = K_t.([0,1//2*t^2]); # rational point

julia> g, _ = transform_to_weierstrass(g, x, y, P);

julia> E6 = elliptic_curve(g, x, y)
Elliptic curve with equation
y^2 + (-t2^2 - 2*t2 - 1)//t2^4*y = x^3
\end{minted}

Let us summarize the equations we found in the following theorem.

\begin{theorem}
Let $Y/\mathbb{C}$ be the singular K3 surface of discriminant 3.
By \cite{Nishiyama1996, FestiVeniani2022} it admits exactly $6$ elliptic fibrations $E_1,\dots E_6$ up to automorphisms.
Weierstrass models of their generic fibers, the Kodaira types of their singular fibers and their Mordell-Weil groups are listed in the following table:\\
\medskip\\
 \begin{tabular}{llll}
 $i$ \quad \phantom{o} & $E_i$ & Kodaira type \hspace{0.5cm} & Mordell-Weil group \\ \hline
 $1$ & $y^2 = x^3 + t^5(t-1)^2$  & $2 \times II^*, IV$ & 0\\
 $2$ & $y^2 = x^3 - t^3x^2 + t^3x $ & $I_6^*$, $III^*$ &  $\ZZ/2\ZZ \times \ZZ$\\
 $3$ & $y^2 = x^3 + (t^3 - 2)x^2 + x $ & $I_{12}^*, I_3$& $\ZZ/2\ZZ$\\
 $4$ & $y^2 = x^3 + t^4x^2 - 8t^2x + 16 $ & $I_{18}$ & $\ZZ/3\ZZ \times \ZZ$\\
 $5$ & $y^2 = x^3 + (t^4 - 8 t)x^2 + 16 t^2 x $& $I_{12}, I_3^*$ & $\ZZ/4\ZZ$\\
 $6$ & $y^2 - (t^2 + 2 t + 1)t^8y = x^3 $ & $3 \times III^*$  & $\ZZ/3\ZZ$
 \end{tabular}
 \medskip
\end{theorem}

To sum up, we have computed all elliptic fibrations up to isomorphism on Vinbergs K3 surface as well as their elliptic parameters and generic fibers in Weierstrass form. Moreover, we can explicitly write down the isomorphisms between the different desingularizations of the Weierstrass models and put them to use to transform between the respective bases of the N\'eron-Severi lattice.

\section*{Acknowledgements}

This work is funded by the Deutsche Forschungsgemeinschaft (DFG, German Research Foundation) – Project-ID 286237555 – TRR 195. The authors wish to thank Matthias Sch\"utt for helpful comments on the manuscript.

\printbibliography

\end{document}